\documentclass[11pt,leqno]{amsart}
\usepackage{amsmath,amsfonts,amssymb,amscd,amsthm,amsbsy}
\newtheorem{thm}{Theorem}[section]
\newtheorem{lem}[thm]{Lemma}
\newtheorem{cor}[thm]{Corollary}

\newtheorem*{THM}{Theorem}

\theoremstyle{definition}
\newtheorem{rmk}[thm]{Remark}
\newtheorem*{defin}{Definition}

\newtheorem*{note}{Note}

\newcommand{\Ht}{\operatorname{ht}}   
\newcommand{\epf}{\operatorname{epf}}  
\newcommand{\rf}{\operatorname{r1f}}   

\author{Raymond C. Heitmann}
\title{The Direct Summand Conjecture in Dimension Three}
\thanks{This reseach was partially supported by a grant from the
  National Science Foundation, DMS-0100731}

\begin{document}  

\maketitle          

In \cite{Ho1}, Hochster proved the following results:

\begin{THM}[Hochster]  
If $R$ is a regular Noetherian ring which contains a field and $S
\supset R$ is a module-finite $R$-algebra, then $R$ is a direct
summand of $S$ as an $R$-module. 
\end{THM}

\begin{THM}[Hochster]  
If $S$ is any local ring which contains a field and $x_1, \ldots, x_n$
is a system of parameters for S, then for every integer $k \geq 0$,
$(x_1 \cdots x_n)^k \not\in \big(x_1^{k+1}, \ldots, x_n^{k+1}\big)S$. 
\end{THM}

The mixed characteristic case of these results is easy for $\dim R
\leq 2$ but relatively little is known for $\dim R>2$.  
The general statements, which are equivalent, became known as the
direct summand and monomial conjectures.  
The principal advance in this subject occurred in Hochster's 1983
article \cite{Ho2} in which he introduced the canonical element
conjecture.  
He demonstrated that this new conjecture was equivalent to the other
two both overall and in the case of fixed dimension.  
He also showed that it was sufficiently strong to imply the validity
of a number of other homological conjectures which had been shown to
follow from the big Cohen-Macaulay modules conjecture. 

In this article, we prove a result (Theorem 3.7) which implies the
three dimensional case of the direct summand conjecture and so also
the three dimensional cases of each of the conjectures that follow
from it.
As it happens, most of these results also follow from the New
Intersection Theorem which was proved by Roberts in the mixed
characteristic case \cite{R1}.
Two exceptions are the improved new intersection conjecture and the
mixed characteristic version of the Evans-Griffith Syzygy Theorem
\cite{EG}.
These two are now theorems in the three dimensional case.
It should also be noted that Hochster \cite{Ho3} has just proved the
existence of balanced big Cohen-Macaulay algebras in dimension three
--- using (3.7) as a key ingredient of the proof.

The main result of this article is of independent interest.
In \cite{H2}, the author introduced several forms of an extended plus
closure. 
These were intended to fill the void due to the absence of a tight
closure analog in mixed characteristic.
A key property of the tight closure is the colon-capturing property.
Theorem 3.7 tells us that the extended plus closures have the
colon-capturing property in dimension three, at least for excellent
rings.
It was shown in \cite{H2} that if any of these new closures had the
colon-capturing property, it would follow that the closures would also
have a second very desirable property --- that all ideals in regular
local rings are closed.
There is an additional implication concerning the Brian\c{c}on-Skoda
Theorem, proved in the general case by Lipman and Sathaye \cite{LS}.
In \cite[p.818]{H2}, mimicking the tight closure treatment, a
generalization of this result was proved for the extended plus
closures.  
Unfortunately, this result is not a true generalization in that it
would only imply the original result if ideals in regular rings were
closed.  
This result now does imply the original in the three dimensional
case. 
So Theorem 3.7 goes a long way towards establishing a good closure
operation for three dimensional mixed characteristic rings.

For the direct summand conjecture, the property that we in fact need
is that ideals in regular rings are closed under the extended plus
closure. This in turn implies that such ideals are closed under the
plus closure.
This latter property is equivalent to the direct summand conjecture.

\section{Preliminaries}

Throughout, $R$ will be a commutative integral domain with unity. 
Also, $R$ will be semi-local of mixed characteristic, that is, the
Jacobson radical of $R$ will contain a prime integer p which is
nonzero as an element of the ring.  
We further assume that each maximal ideal of $R$ has the same height.  
$R^+$ will denote the integral closure of $R$ in an algebraic closure
of its quotient field.  
We will refer to $x_1,\ldots, x_n$  as a set of parameters in $R$
provided $\Ht(x_{j_1}, \ldots x_{j_k})R=k$ for any $k$ distinct
elements of the set.   
If $R$ is catenarian, this is the same as assuming that
$\Ht(x_1,\ldots, x_n)R =n$ or that the set is part of a complete
system of parameters. 
By $H_1(x_1, x_2, x_3; R)$, we mean the usual Koszul homology.  
In order to apply our main theorem, we need to recall a critical
result from \cite[p.61]{HH}.

\begin{thm} \label{thm:1.1}
If $R$ is an excellent biequidimensional semi-local ring, $c \in R$
and $R_c$ is
Cohen--Macaulay, then $c$ has a power $c^k$ such that for every set of
parameters $x_1,\ldots,x_n$, \break $c^kH_i(x_1,\ldots,x_n;R)=0 $ for
all $i>0$. 
\end{thm}

In particular, when $R$ is integrally closed of dimension three, we know
that $R_c$ is Cohen-Macaulay for every non-unit $c$ and so every element 
of the maximal ideal has a power which kills all of the first Koszul
homology modules.

Next we recall some definitions from \cite{H2}.  
The last two were introduced in that article; the first has a longer
history.  
Since we are only concerned with integral domains which do not contain
the rational numbers, we may state the definitions more simply.

\begin{defin}
If $x \in R$, then $x$ is in the plus closure of $I$ if $x \in IR^+
\cap R$.  
We write $x \in I^+$.
\end{defin}

\begin{defin}
If $x \in R$, then $x$ is in the full extended plus closure of $I$ if
there exists $c \neq 0 \in R$ such that for every positive integer
$n$, $c^{1/n}x \in (I, p^n)R^+$.  
We write $x \in I^{\epf}$.
\end{defin}

\begin{defin}  
If $x \in R$, then $x$ is in the full rank one closure of $I$ if for
every rank one valuation on $R^+$, every positive integer $n$, and
every $\varepsilon > 0$, there exists  $d \in R^+$ with $v(d) <
\varepsilon$  such that $dx \in (I, p^n)R^+$.  
We write $x \in I^{\rf}$.
\end{defin}

In this article, we shall focus almost exclusively on the full
extended plus closure.  
In fact, we could restrict our attention to a smaller closure
which requires only that $c^{1/n}x \in IR^+$.  
However, this smaller closure is certain to be too small to serve as a
tight closure analog in dimension four or larger.  
Unfortunately, in focusing on dimension three, this article will shed
no light on whether or not the $p^n$ is sufficiently helpful.  
We should note the trivial relationship between the closures --- $I^+
\subset I^{\epf} \subset I^{\rf}$. 

Next we cite the critical result from \cite[p.820]{H2}.

\begin{thm} \label{thm:1.2}
Let $R$ be a regular local ring with $p \in J(R)$.  
If either full rank one closure or full extended plus closure has the
colon-capturing property for finite extensions of regular local rings,
then  $I^{\epf} = I^{\rf} =I$ for every ideal $I$ of
$R$. 
 In particular, $I^+ = I$ for every ideal $I$ of $R$.
\end{thm}

The statement in the earlier article did not mention the full extended
plus closure in the hypothesis.  
However, it does not change the theorem since $I^{\epf} \subset
I^{\rf}$ and it emphasizes how we shall use it here.  
There are two strengthenings of this result which we shall need and
which are readily observable from reading the proof.  
The first is the very logical observation that in order to prove that
ideals are closed in a particular regular local ring, we need only
have the colon-capturing property for finite extensions of that ring.  
Hence the theorem may be restated for excellent rings of a fixed
dimension.  
Moreover, the final statement still holds for all regular local rings
of that dimension since if we were to have an ideal which was not
plus-closed in a regular local ring, we would have such an ideal in
its completion.  
The second observation, slightly more subtle, is that we only needed
the colon-capturing property for systems of parameters which are
powers of a single system of parameters.  
Rather than relying on this observation, we shall prove that, in the
primary situation of interest here, this is the same as
colon--capturing in general. 

\begin{lem} \label{lem:1.3}
Suppose $\dim R=3$ and   $y_1, y_2, y_3$ is a system of parameters.  
Further suppose that for every system of parameters $x_1, x_2, x_3$,
each $y_i$ kills $H_1(x_1, x_2, x_3;R)$.  
If for some $\delta \in R^+$, we have $\delta\big((y_1, y_2):_R y_3\big) \subset
(y_1, y_2)R^+$, then $\delta\big((x_1, x_2):_R x_3\big) \subset(x_1,
x_2)R^+$. 
\end{lem}

\begin{proof} 
In this proof we shall use Roman letters for elements of $R$ and Greek
letters for elements in $R^+$.
First we choose an element $r$ such that $z_2=y_2+ry_1$ is not
contained in any height one prime ideal which contains $x_1$. 
Next choose elements $r_1,r_2$ such that $z_3=y_3+r_1y_1+r_2y_2$ is
contained in no height two ideals which contain  $(x_1,x_2)R$ or
$(x_1,z_2)R$.   
So  $x_1,x_2,z_3$ and $x_1, z_2, z_3$ are also systems of parameters
as is $y_1, z_2, z_3$. 

Suppose $a_3 \in \big((x_1,x_2):_Rx_3\big)$.  
Then $a_3 \in \big((x_1,x_2):_Rz_3\big)$ since $z_3$ kills colons and
so we may write $a_1 x_1 + a_2 x_2 + a_3 z_3 =0$.  
As $a_2 \in \big((x_1,z_3):_Rx_2\big) \subset \big((x_1,z_3):_Rz_2\big)$,
we may write $b_1 x_1 + a_2 z_2 + b_3 z_3 = 0$.  
Again, as $b_1 \in \big((z_2, z_3):_R x_1\big) \subset \big((z_2, z_3):_R
y_1\big)$, we may write $b_1 y_1 + c_2 z_2 + c_3 z_3 = 0$.  
Since $(y_1, z_2)R = (y_1, y_2)R$, we have $\delta c_3 \in (y_1, z_2)R^+$.  
Writing $\delta c_3 = \alpha y_1 + \beta z_2$, we have a relation
$(\delta b_1+\alpha z_3)y_1 + (\delta c_2 + \beta z_3)z_2 = 0$ and since no
height one prime ideal of $R^+$ contains $y_1$ and $z_2$, $\delta b_1 + \alpha z_3
\in z_2 R^+$; hence $\delta b_1 \in (z_2, z_3)R^+$.   
In an identical fashion, we see that $\delta a_2 \in (x_1, z_3)R^+$ and
finally that $\delta a_3 \in (x_1, x_2)R^+$.
\end{proof}

\begin{cor} \label{cor:1.4}
Suppose $\dim R=3$, and $y_1, y_2, y_3$
is a system of parameters.  
Further suppose that for any system of parameters $x_1, x_2, x_3$,
each $y_i$ kills $H_1(x_1, x_2, x_3;R)$.  
If for some $c \in R$, we have  $c^{1/n} \big((y_1, y_2):_R y_3\big)
\subset (y_1, y_2)R^+$ for every $n$, then the full extended plus
closure has the colon-capturing property for $R$.
\end{cor}

\begin{proof} 
Clearly, even the potentially smaller plus closure captures colons
of the form $(0:x)$ and $(x_1: x_2)$.  
The hard case is given by the lemma.
\end{proof}

Next we shall define a function $\tau$ on the positive integers which
shall play a critical role in the proof of the main theorem.  
In light of (\ref{lem:1.6}), it seems likely that this function has
been used previously; however we have not seen it.  
Strictly speaking, $\tau$ depends on the choice of prime integer $p$
and should truly be denoted $\tau_p$, but as the function will only be
employed for $p$ the characteristic of the residue field, we shall
omit the subscript.  

\begin{defin}  
For a positive integer $n$, express $n-1$ in base $p$ and let
$\bar{\tau}(n)$ be the sum of the digits.  
We take $\bar{\tau}(1)=0$.  
Then define 
$$
\tau(n) = \frac{\bar{\tau}(n)}{p-1}\ .
$$
\end{defin}

\begin{note}
It is more natural to define $\tau$ without the shift --- say define
$\sigma$ by $\sigma(n-1) = \tau(n)$ --- and no doubt those who have
encountered this function previously have seen it in this form.
In fact, the statement of (1.6) is nicer using $\sigma$.
However, the proof of the main theorem is simpler using $\tau$ and in
fact the shifted function appears to better capture what is going on.
\end{note}

\begin{lem} \label{lem:1.5}
Suppose $i,j$ are positive integers and $n=i+j$.  
Let  $a_k \cdots a_0$ be the expression for $i$ in base $p$, i.e.,
$i= a_0 + a_1 p + \cdots + a_kp^k$  with $0 \leq a_j < p$.  
Similarly, suppose $j = b_k \cdots b_0$, $n =c_k \cdots c_0$.  
Let $d= \big\lvert \{ j \mid a_j + b_j > c_j\} \big\rvert$.  
Then the highest power of $p$ which divides $\binom{n}{i}$ is $p^d$. 
\end{lem}

\begin{proof} 
This is just a restatement of a lemma in \cite[p.696]{H1}.  
It is actually slightly stronger than the statement there but is in
fact what is actually proved.
\end{proof}

\begin{lem} \label{lem:1.6}  
Let  $0 < j < i < p^L$ be integers.
\begin{itemize}
\smallskip
\item[(a)] The highest power of $p$ which divides $\binom{i-1}{j-1}$ is
  $\tau(j) + \tau(i-j+1) - \tau(i)$.
\smallskip
\item[(b)] The highest power of $p$ which divides $\binom{p^L-j}{i-j}$
  is also $\tau(j) + \tau(i-j+1) - \tau(i)$.
\smallskip
\item[(c)] The highest power of $p$ which divides $\binom{p^L}{i}$ is
  $L + \tau(i+1) - \tau(i) - \tau(2)$.
\end{itemize}
\end{lem}

\begin{proof} (a)  To perform the addition $(i-1)+(n-i)=n-1$ in base $p$, one can 
first add the corresponding digits, allowing digits as large as
$2p-2$.  
The sum of the digits of this ``number'' will clearly be
$\bar{\tau}(j) + \bar{\tau}(i-j+1)$.  
Next, starting from the right, for each digit which is at least $p$,
we deduct $p$ and add one to the next digit.  
Each such carrying operation clearly reduces the sum of the digits by
$p-1$.  
Thus the highest power of $p$ which divides  $\binom{i-1}{j-1}$ times
$(p-1)$ is $\bar{\tau}(j) + \bar{\tau}(i-j+1) - \bar{\tau}(i)$.  
The result quickly follows.

(b) This will follow from (a) once we verify that the highest power
of $p$ which divides  $\binom{p^L-j}{i-j}$ is the same as the highest
power of $p$ which divides $\binom{i-1}{j-1}$.  
To see this, note that $\binom{i-1}{j-1} = \binom{i-1}{i-j}$ and the
two quantities are $\frac{(p^L-j) \cdots (p^L-(i-1))}{(i-j)!}$
and $\frac{(i-1)\cdots j}{(i-j)!}$ respectively. \smallskip 
Finally, whenever $j \leq k \leq i-1$, $k<p^L$, and a power of $p$
will divide $p^L-k$ if and only if it divides $k$. 
The result is now clear.

(c) We note that the highest power of p which divides $\binom{p^L}{i}$
is the highest power of $p$ which divides $\frac{p^L}{i} =
\frac{p^L}{\binom{i}{1}}$. 
The result now follows from
part (a).
\end{proof}

We conclude this section with two lemmas that will be needed for
the main proof.

\begin{lem} \label{lem:1.7}  
Let $K \geq -1$ be an integer and let $e = \frac{1}{p^{K+1}}$. 
Then $ei + K - \tau(i) \geq 0$ for every positive integer $i$.
\end{lem}

\begin{proof}
We induct on $K$.
For $K=-1$, $e=1$ and so we must show $i - 1 - \tau(i) \geq 0$ or
equivalently $\tau(i) \leq i -1$.
But $\tau(i)$ is simply the sum of the digits of $i-1$ divided by a
positive integer and so cannot possible exceed $i-1$.
Now we consider $K>-1$ and assume the result holds for $K-1$.
Suppose there exists $i$ with $ei + K -\tau(i) < 0$.
If $p$ does not divide $i$, then $\tau(i+1) = \tau(i) + \frac{1}{p-1}$
while $e(i+1) = ei + e  \leq ei + \frac1p$.
It follows that $e(i+1) + K - \tau(i+1) < 0$.
Similarly, if $p$ does not divide $i+1$, we see that $i+2$ provides a
counterexample.
Thus we may produce a counterexample to the lemma with $p$ dividing $i$.
Writing $i=pj$ we have $epj + K - \tau(pj) < 0$.
Since the base $p$ representation of $pj-1$ is just the base $p$
representation of $j-1$ with the final digit $p-1$ added, we have
$\tau(pj) = \tau(j) + 1$.
Thus we have $\frac{p}{p^{K+1}}j + K - 1 - \tau(j) < 0$.
However, this contradicts the $K-1$ case of the lemma.
\end{proof}

\begin{lem} \label{lem:1.8}  
Let $R$ be an integrally closed domain which satisfies our usual
assumptions and suppose that there is an element $\sigma \in R$ such
that $\sigma^{p-1} = p$.
Let $x,y,z \in R$ where no height one prime ideal contains both $p$
and $y$.
Suppose $z = xw + yv$ where $w,v$ are integral over $R[p^{-1}]$ and
$w$ satisfies the monic polynomial $T^n + a_1 T^{n-1} + \cdots + a_n$.
Further suppose $a_i \in p^{K- \tau(i)}R$ for every $i$.
If $e = \frac1{p^{K+1}}$, then $p^ew$ and $p^ev$ are integral over
$R$.
\end{lem}

\begin{proof}
If necessary, we adjoin a $(p^{K+1})$st root of $p$ and denote it
$p^e$.
Then $p^ew$ satisfies the monic polynomial $T^n + d_1 T^{n-1} + \cdots
+ d_n$ with $d_i = p^{ei} a_i$ for every $i$.
Thus $d_i \in p^{ei + K - \tau(i)}R$.
By (1.7), $d_i \in R$ and so $p^ew$ is integral over $R$.
Now this implies that $p^ev$ is integral over $R[y^{-1}]$.
By hypothesis, $p^ev$ is also integral over $R[p^{-1}]$.
Since $R$ is integrally closed and no height one prime ideal contains
both $p$ and $y$, we see that $p^ev$ is integral over $R$. 
\end{proof}

\section{Main Results}
In the succession of lemmas leading up to the main theorem, $S$ will
denote an integral domain which contains the rational numbers.
The results remain true and the proofs remain valid if $S$ contains
$\mathbb{Z}_{(p)}$, but in fact all applications of the lemmas will
be to the ring $S=R[p^{-1}]$.

Demonstrating colon-capturing in dimension three amounts to proving 
that certain elements are in the plus closure of a two-generated ideal. 
A key tool for doing this is the following lemma from \cite[p.693]{H1}.

\begin{lem} \label{lem:2.1}  
Suppose $x,y,z\in S$ with $y\ne 0$. 
Let 
$$
f(T) = \sum_{i=0}^n a_i T^{n-i}
$$ 
be a monic polynomial over $S$ and suppose $w$ is an element in an 
extension domain of $S$ such that ${f(w)=0}$.  
For $0 \leq i \leq n$, set  
$$
b_i = (-1)^i \sum_{j=0}^i \binom{n-j}{i-j} a_j z^{i-j}x^j
$$
and let 
$$
g(T)=\sum_{i=0}^n b_i T^{n-i}\ .
$$  
Then $g(z-xw)=0$.  
In particular, if each $b_i \in y^iS$, $(z-xw)/y$ is integral over
$S$. 
\end{lem}

In order to prove the main theorem, we have to get from our starting
point, an equation $p^Nz = cx + dy$ in $R$, to an equation $\rho z =
\gamma x + \delta y$ where $\gamma, \delta$ are integral over $R$ and
$\rho$ is an $m$th root of $p$ for some large $m$.
Equivalently, we need to write  $z=p^{-N}cx + p^{-N}dy$ as
$\rho^{-1}\gamma x + \rho ^{-1}\delta y$.
The element $\rho^{-1}\gamma$ will be the $w$ of (2.1) and $w$ will be
defined indirectly, by constructing $f(T)$ and letting $w$ be a root
of $f(T)$.
According to the basic plan (which will have to be modified somewhat),
the polynomial $f(T)$, with coefficients in $R[p^{-1}]$, will be
constructed via a recursive procedure, which successively defines
$a_1, a_2, \ldots, a_n$.
By choosing $a_1$ such that $b_1 \in yR[p^{-1}]$, $a_2$ such that $b_2
\in y^2 R[p^{-1}]$, etc., we will force $\rho^{-1} \delta$ to be
integral over $R[p^{-1}]$.

If we omit the condition that $\rho w$ be integral over $R$,
constructing $f(T)$ is easy.
In fact $(T - p^{-N}c)^n$ works for any choice of $n$.
The next lemma is critical for choosing a better $f(T)$.
It asserts that there is no subtle way that the recursion can fail.
No matter how we select an initial sequence of $a_i$'s (subject of
course to the conditions on the $b_i$'s), we can always finish it if
we omit the condition that $\rho w$ be integral over $R$.

\begin{lem} \label{lem:2.2}
Suppose $z=ax+by$ with $a,b,x,y \in S$ and let $n$ be a positive
integer.
Further suppose $a_0=1, a_1,
\ldots, a_{k-1} \in S$ have been chosen with $k \leq n$ so that 
$$
\sum_{j=0}^i \binom{n-j}{i-j} a_j z^{i-j} x^j \ \in \ y^i S
\quad \text{for} \quad i=1, \ldots, k-1\ .
$$ 
Then we may find $a_k \in S$ such that 
$$
\sum_{j=0}^k \binom{n-j}{k-j} a_j z^{k-j} x^j \ \in \ y^k S\
. 
$$.
\end{lem}

\begin{proof} 
First we will show that we can reduce to the case $z=ax$.
As the lemma can be applied repeatedly if true, it is equivalent to
the assertion that we may finish the polynomial $f(T) =
\sum_{i=0}^{k-1}a_i T^{n-i}$
to obtain a polynomial  $f(T) =  \sum_{i=0}^n  a_i T^{n-i}$ 
such that for every root $w$ of $f(T)$, $(z-xw)/y$ is integral over
$S$.   
Since $z=ax+by$, this is the same as $(ax+by-xw)/y$ being integral
over $S$ and of course this is that same as $(a-w)x/y$ being integral
over $S$.
Thus we can make the reduction.
Now, returning to the original statement of the conclusion of the
lemma, we need only be able to find $a_k$ such that 
$$
\sum_{j=0}^k \binom{n-j}{k-j}a_j a^{k-j} x^k\in y^kS
$$
or equivalently such that
$$
\sum_{j=0}^k \binom{n-j}{k-j}a_j a^{k-j} \ \in \ (y^k:_Sx^k)\ .
$$  
This is trivial since $a_k$ occurs on the left hand side with
unit coefficient.
\end{proof}

We are now ready to prove a special case of the main result.
The proof of this proposition illustrates several of the ideas which
form the basis of the later result.

\begin{thm}
Let $p,x,y$ be parameters in $R$, an integrally closed domain
which satisfies our usual assumptions.  
Assume there is an element $\sigma \in R$ with $\sigma^{p-1}=p$ and
that $p^N,x,y$ kill $H_1(p^m, x^m, y^m; R)$ for every positive integer
$m$.   
Further assume that $H_1(p^N, x, y; R)$ is cyclic, generated by
$(z,a,b)$.
Then for any rational $e>0$, there is a module-finite extension $A$
of $R$ with $p^ez\in(x,y)A$.
Thus $z \in (x,y)R^{\epf}$.
\end{thm}

\begin{proof}  
It suffices to prove the result with $e = \frac{1}{p^{K+1}}$ for
arbitrary $K$. 
We set $L=K+N$.
We shall construct a polynomial $f(T) =  T^{p^L} + a_1 T^{p^L-1} +
\cdots + a_{p^L}$ with coefficients in $R[p^{-1}]$ such that if $w$ is
any root of $f(T)$, $v = (z-xw)/y$ is also integral over $R[p^{-1}]$.
If we can accomplish this with each $a_j \in p^{K-\tau(j)}R$, then by
(1.8), the conclusion holds with $A = R[p^e, p^ew, p^ev]$.
Thus the entire proof rests upon our ability to satisfactorily choose
the $a_j$'s.

We shall eventually choose elements $a_i', a_i'' \in p^{K-\tau(i)}R$
for every integer $i$, $1 \leq i \leq p^L$.
We shall define $a_0 = 1$ and $a_i = a_i'+y^ia_i''$ for every $i \geq 1$.
Then, by (2.1), $v = (z-xw)/y$ is integral over $R[p^{-1}]$ provided
$\sum_{j=0}^i \binom{p^L-j}{i-j} a_j z^{i-j}x^j \in y^i R[p^{-1}]$ for
each $i\leq p^L$.
We shall refer to this inclusion for a specific $i$ as the $i$th
condition.
The elements will be chosen to satisfy the conditions by a recursive
procedure beginning with $i=1$.
At each step, we choose $a_i' \in  p^{K-\tau(i)} R$ and $a_{i-1}'' \in
p^{K-\tau(i-1)} R$ such that the $i$th condition is satisfied.
The $i=1$ step is irregular as there is no $a_0''$.
Finally, the superfluous $a_{p^L}''$ equals zero.

Since $p^Nz +ax +by =0$, $p^{N+K}z + p^Kax \in yR[p^{-1}]$.
Hence, letting $a_1'=p^Ka \in p^K R = p^{K-\tau(1)}R$, we have $p^Lz +
a_1 'x \in yR[p^{-1}]$ and the first condition is satisfied regardless
of the choice of $a_1''$.
In fact, it is important to note that the choice of $a_i''$ never
affects the $i$th condition.
For the recursive step, we assume that we have chosen each $a_k$
properly for $k<i-1$ as well as $a_{i-1}'$ so that the first $i-1$
conditions are satisfied.
For notational simplicity, it is convenient to assume $a_{i-1}
=a_{i-1}'$.
As we know that  $\sum_{j=0}^k \binom{p^L-j}{k-j} a_j z^{k-j}x^j \in
y^k R[p^{-1}]$ for each $k<i$, (2.2) yields
$$
\sum_{j=0}^{i-1} \binom{p^L-j}{k-j} a_j z^{k-j}x^j \in (x^i,y^i)
R[p^{-1}]\ . 
$$   
Now let $z_i = p^{-K + \tau(i) -\tau(2)}\sum_{j=0}^{i-1}
\binom{p^L-j}{i-j} a_j z^{i-j}x^j$.
We claim $z_i \in R$.
To prove the claim, it suffices to prove it one term at a time and so
we need only show $p^{-K + \tau(i) -\tau(2)}\binom{p^L-j}{i-j} a_j \in
R$ for each $j < i$.
Using (1.6) and the fact that $a_j \in p^{K-\tau(j)}R$, we see that
this is equivalent to 
$$
\big(-K + \tau(i) - \tau(2)\big) + \big(\tau(j) + \tau(i-j+1)
-\tau(i)\big) + \big(K-\tau(j)\big) \geq 0
$$
when $j>0$.
But this is clear since the left hand side is just $\tau(i-j+1) -
\tau(2)$.
For $j=0$, we need
$$
\big(-K + \tau(i) - \tau(2)\big) + \big(L + \tau(i+1)
-\tau(i) - \tau(2)\big) \geq 0\ .
$$
Since $L \geq K+1$, it suffices to show that $1 + \tau(i+1) \geq 2
\tau(2)$ and this too is clear.
So the claim is verified.

As $z_i \in (x^i,y^i)R[p^{-1}] \cap R$, it follows that $z_i \in
\big((x^i,y^i):_Rp^m\big)$ for some $m$.
Since $x,y,p^N$ kill $H_1(p^m,x^m,y^m;R)$ for every positive integer
$m$, $z_i \in (x^i,y^i)R + (xy)^{i-1}\big((x,y):_R p^N\big)$.
Then, by the cyclicity assumption, $z_i \in (x^i, y^i, (xy)^{i-1}z)R$.
Write $z_i = c_{i1}x^i + c_{i2}y^i + c_{i3}(xy)^{i-1}z$.
Since $\sum_{j=0}^{i-1} \binom{p^L-j}{i-j} a_j z^{i-j}x^j =
p^{K-\tau(i)+ \tau(2)}z_i$, the choices 
$$
a_{i-1}'' = -\frac{p^{K-\tau(i)+ \tau(2)}c_{i3}}{p^L-(i-1)} \qquad
\text{ and } \qquad a_i' = -p^{K-\tau(i)+ \tau(2)}c_{i1}
$$
satisfy the $i$th condition.
Obviously, $a_i' \in p^{K-\tau(i)}R$.
Finally, to see that $a_{i-1}''\in p^{K-\tau(i-1)}R$, we note that the
highest power of $p$ which divides $p^L - (i-1)$ is the highest power
dividing $\binom{i-1}{1}$, namely $p^{\tau(i-1)+ \tau(2) - \tau(i)}$.
Thus it suffices to show that $\big(K-\tau(i)+ \tau(2)\big) -
\big(\tau(i-1) + \tau(2) - \tau(i)\big) \geq K - \tau(i-1)$. 
This is an equality and the recursion step is complete, finishing the
proof. 
\end{proof}

The cyclicity assumption only came into play at one point --- to give
us $z_i \in (x^i, y^i, (xy)^{i-1}z)R$.
The key to dropping that assumption is to develop a way of using
earlier $z_i$'s to help with later ones.
We need to be able to say that for sufficiently large $n$ and
$i>n$, $z_i \in (x^i, y^i, (xy)^{i-1}z_1, (xy)^{i-2}z_2, \ldots,
(xy)^{i-n}z_n)R$. 
This will be shown presently; however a number of complications lay
ahead.
Instead of choosing each $a_i$ in a two step procedure with an initial
choice and then a single adjustment to help at the $(i+1)$st step, we
must allow for multiple adjustments --- adjustments that don't
obviously preserve previously satisfied conditions.
It is also not clear at the start of the process what the correct
degree of the polynomial should be and so we must be able to adjust
this as we go along.

\begin{lem} 
Let $p,x,y$ be parameters in $R$, an integrally closed domain which
satisfies our usual assumptions.
Assume that $x,y$ kill $H_1(p^m,x^m,y^m;R)$ for every positive integer
$m$.
Further assume that $z_i \in \big((x^i,y^i):_Rp^N\big)$ for every
positive integer $i$.
Then there exists a positive integer $n$ such that for all $i > n$, 
$z_i \in \big(x^i, y^i, (xy)^{i-1}z_1, (xy)^{i-2}z_2, \ldots,
(xy)^{i-n}z_n \big) R$.
\end{lem}

\begin{proof}
Let $M_i = \big((x^i,y^i):p^N\big)$. 
Since $(x,y)\big((x^i,y^i):p^N\big) \subset (x^i,y^i)R$ and $R$ is
integrally closed, it follows that $\big((x^i,y^i):p^N\big) \subseteq
(x^i,y^{i-1})R  \cap (x^{i-1},y^i)R = (x^i, y^i, x^{i-1}y^{i-1})R$.
So $M_i = x^iR + y^iR + (xy)^{i-1}M_1$.
Also let $Q_i = \big(x^i, y^i, (xy)^{i-1}z_1, (xy)^{i-2}z_2, \ldots,
z_i\big) R$.
Clearly $Q_i \subset M_i$ and, noting $xyM_{i-1} \subset M_i$,
$xyQ_{i-1} \subset Q_i$, multiplication by $xy$ induces the commutative
diagram
\begin{equation*}
\begin{CD}
Q_i & \quad\subset\quad & M_i \\
@AAA   @AAA \\
Q_{i-1} & \quad\subset\quad & M_{i-1} 
\end{CD}\ .
\end{equation*}
Next Let $\overline{M}_i$ and $\overline{Q}_i$ be the respective
quotients modulo $(x^i, y^i)R$.
Then our previous diagram induces
\begin{equation*}
\begin{CD}
\overline{Q}_i & \quad\subset\quad & \overline{M}_i \\
@AAA   @AAA \\
\overline{Q}_{i-1} & \quad\subset\quad & \overline{M}_{i-1} 
\end{CD}\ .
\end{equation*}
But the map on the right is an isomorphism and so $\overline{Q}_1 \to
\overline{Q}_2 \to \cdots \to \overline{Q}_i \to$ may be viewed as an
ascending chain of submodules of the same noetherian module.
Thus for some $n$, $\overline{Q}_{i-1} \to \overline{Q}_i$ is
surjective for all $i>n$ and the result follows.
\end{proof}

Before going on to the main theorem, we need several technical lemmas
and the ideas behind them.
First we remark that the easiest proof of (2.1) is a simple
application of Taylor's Theorem.
In fact, $b_i$ is just $\frac{(-1)^ix^i}{(n-i)!}
f^{(n-i)}\big(\frac{z}{x}\big)$ where $n$ is the degree of $f(T)$.
As $S$ contains the rationals, the condition $b_i \in y^iS$ is
equivalent to $x^if^{(n-i)}\big(\frac{z}{x}\big) \in y^i S$.
Interestingly, for fixed $i$, the latter condition is stable under
taking antiderivatives.
To see this, note that if $\tilde{f}(T)$ is an antiderivative of
$f(T)$, then the degree of $\tilde{f}(T)$ is $n+1$ and
$x^i\tilde{f}^{(n+1-i)} \big(\frac{z}{x}\big) = x^i f^{(n-i)}
\big(\frac{z}{x}\big)$.
Likewise, the condition is stable under taking derivatives provided
$i<n$.

In the proof of the main theorem, we will have occasion to exploit
both these observations, the first to increase the degree of the
polynomial without destroying the previous steps and the second to
``improve'' previously chosen coefficients.
The next two technical lemmas are just combinatorial translations of
the two stability conditions, stated in a manner so that they can be most
easily applied here.
Lemma 2.5 relates to the antiderivative while the $M=L$ case of (2.6)
relates to the derivative.
The full statement of (2.6) is a combination.

\begin{lem}
Suppose $M>L$ are positive integers and $x,y,z, a_1, \ldots, a_k$ are
elements of $S$ such that $\sum_{j=0}^i \binom{p^L-j}{i-j} a_j
z^{i-j}x^j \in y^i S$ for all $i \leq k$ where $k$ is an integer less
than $p^L$.
If $\tilde{a}_j = \big(\prod_{m=0}^{j-1}\frac{p^M-m}{p^L-m}\big)a_j$
for each $j>0$, then $\tilde{a}_j = p^{M-L}q_ja_j$ where $q_j$ is a
unit in $\mathbb{Z}_{(p)}$ and $\sum_{j=0}^i \binom{p^M-j}{i-j}
\tilde{a}_j z^{i-j}x^{j} = p^{M-L}q_i \sum_{j=0}^i \binom{p^L-j}{i-j}
a_j z^{i-j}x^{j} \in y^iR$ for all $i\leq k$.
\end{lem}

\begin{proof}
If $0 \neq m < p^L$, the highest power of $p$ that divides $p^L-m$ is
the highest power of $p$ that divides $m$ and the same is true for
$p^M-m$ and so $\frac{p^M-m}{p^L-m}$ is a unit in $\mathbb{Z}_{(p)}$.
For $m=0$, the quotient is of course $p^{M-L}$.
This gives the first half of the conclusion.
Next note that 
$$
\frac{\binom{p^M-j}{i-j}}{\binom{p^L-j}{i-j}} = 
\frac{(p^M-j) \cdots (p^M-i+1)}{(p^L-j) \cdots (p^L-i+1)}\ .
$$
Hence
$$
\binom{p^M-j}{i-j}\tilde{a}_j = \Big(\prod_{m=0}^{i-1}
\frac{p^M-m}{p^L-m}\Big) \binom{p^L-j}{i-j}a_j.
$$
The second half follows immediately.
\end{proof}

\begin{lem}
Let $d,i,M,L$ be integers with $L \leq M$ and $0 < d < i \leq p^L$.
Suppose $x,y,z, a_0, a_1, \ldots, a_{i-1}$ are elements of $S$ such
that 
$$
\sum_{m=0}^h \binom{p^L-m}{h-m} a_m z^{h-m}x^m \in y^h S
$$
for all integers $h < i$.
Let $\tilde{a}_j=0$ for $j<d$ and 
$$
\tilde{a}_j = \frac{y^d \binom{p^L-j+d}{i-j} a_{j-d}}{\binom{p^M-j}{i-j}}
$$
for $d \leq j \leq i-1$.
Then
$$
\sum_{j=0}^k \binom{p^M-j}{k-j}\tilde{a}_j z^{k-j}x^j \in y^kS
$$
for all $k< i$.
\end{lem}

\begin{proof}
For every $k<i$, we have
\begin{equation*}
\begin{split}
\sum_{j=0}^k \binom{p^M-j}{k-j}\tilde{a}_j z^{k-j}x^j & =
\sum_{j=d}^k \binom{p^M-j}{k-j}y^d a_{j-d}\binom{p^L-j+d}{i-j}
z^{k-j}x^j / \binom{p^M-j}{i-j} \\
  & = y^d \sum_{m=0}^{k-d} \binom{p^M-m-d}{k-m-d}a_m
  \binom{p^L-m}{i-m-d} z^{k-d-m}x^{m+d} / \binom{p^M-m-d}{i-m-d}\ .
\end{split}
\end{equation*}
Next we compute
\begin{equation*}
\begin{split}
\frac{\binom{p^M-m-d}{k-m-d}
  \binom{p^L-m}{i-m-d}}{\binom{p^M-m-d}{i-m-d}} & =
\frac{(p^M-m-d)!(p^L-m)!(i-m-d)!(p^M-i)!}
  {(k-m-d)!(p^M-k)!(i-m-d)!(p^L-i+d)!(p^M-m-d)!} \\
  & = \frac{(p^L-m)!(p^M-i)!}{(k-m-d)!(p^M-k)!(p^L-i+d)!} \\
  & = \binom{p^L-m}{k-d-m}q  
\end{split}
\end{equation*}
where $q$ is a rational number independent of $m$.
Thus
$$
\sum_{j=0}^k\binom{p^M-j}{k-j}\tilde{a}_j z^{k-j}x^j  =
y^dx^d q \sum_{m=0}^{k-d} \binom{p^L-m}{k-d-m}a_m z^{k-d-m}x^m\ .
$$
Since the final sum is in $y^{k-d}S$ by hypothesis, the first sum is
in $y^kS$ as desired.
\end{proof}

We are at last ready to prove the main result.

\begin{thm} 
Let $p,x,y$ be parameters in $R$, an integrally closed domain
which satisfies our usual assumptions.  
Assume there is an element $\sigma \in R$ with $\sigma^{p-1}=p$ and
that $p^N,x,y$ kill $H_1(p^m, x^m, y^m; R)$ for every positive integer
$m$.   
Suppose $p^Nz \in (x,y)R$.  
Then for any rational $e>0$, there is a module-finite
extension $S$ of $R$ with $p^ez\in(x,y)S$.  
Thus $z \in (x,y)R^{\epf}$.
\end{thm}

\begin{proof}
The basic structure of the proof echoes that of (2.3).
We prove the result with $e=\frac{1}{p^{K+1}}$ for arbitrary $K$ by
constructing a polynomial $f(T) = T^{p^L} + a_1 T^{p^L-1} + \cdots +
a_{p^L}$  with coefficients in $R[p^{-1}]$ such that if $w$ is any
root of $f(T)$, $v=(z-xw)/y$ is also integral over $R[p^{-1}]$.
If we can accomplish this with each $a_j \in p^{K - \tau(j)} R$, then
by (1.8), the conclusion holds with $S=R[p^e,p^ew,p^ev]$.
Thus the entire proof rests on our ability to satisfactorily choose
the $a_i$'s.
However, unlike the (2.3) proof, $L$ is not determined at the start;
it will also be chosen in the recursive procedure.

We now describe the recursion.
For each integer $i$, $1 \leq i \leq p^L$, we shall choose a positive
integer $L_i \geq L_{i-1}$ and elements $a_{1i}, \ldots, a_{ii}$ such
that (1) $a_{ji} \in p^{K-\tau(j)} R$ for every $j$ and (2), with
$a_{0i}=1$, the condition
$$
\sum_{j=0}^k \binom{p^{L_i}-j}{k-j}a_{ji}z^{k-j}x^j \in y^kR[p^{-1}]
$$
will be satisfied for each $k \leq i$.
The procedure ends when $i = p^{L_i}$ (something we still must
demonstrate happens), at which time we get the desired polynomial
$f(T)$ with $L=L_i$ and $a_j = a_{ji}$ for every $j$.
We also let $z_1 = z$ and for $i >1$, we let 
$$
z_i = p^{-E_i} \sum_{j=0}^{i-1} \binom{p^{L_{i-1}}-j}{i-j} a_{j, i-1}
z^{i-j} x^j
$$
where $E_i = K - \tau(i) + \tau(2)$.
We shall also show that (3) $z_i \in \big((x^i,y^i):_R p^N\big)$.

For the initial step $(i=1)$, we let $D_1$ be the smallest integer
such that $p^{D_1}z \in (x,y)R$.
Then we may find $a \in R$ such that $p^{D_1}z + ax \in yR$.
Choose $L_1 = K+D_1$ and $a_{11} = p^Ka$ and we see that the first
condition is satisfied with $a_{11} \in p^KR = p^{K-\tau(1)}R$.
By hypothesis, $z_1 \in  \big((x,y):_R p^N\big)$.

For $i>1$, we first demonstrate (3).
This is identical to the situation in (2.3).
As we know that
$$
\sum_{j=0}^k \binom{p^{L_{i-1}}-j}{k-j} a_{j, i-1} z^{k-j} x^j \in
y^kR[p^{-1}] 
$$
for each $k<i$, (2.2) yields
$$
\sum_{j=0}^{i-1}\binom{p^{L_{i-1}}-j}{k-j} a_{j, i-1} z^{k-j} x^j \in
(x^i,y^i)R[p^{-1}]  
$$
and so $z_i \in (x^i,y^i)R[p^{-1}]$.
To see that $z_i \in R$, it suffices to prove it one term at at time
and so we need only show $p^{-K+\tau(i)-\tau(2)}
\binom{p^{L_{i-1}}-j}{i-j}a_{j,i-1} \in R$ for each $j<i$.
Using (1.6) and the fact that $a_{j,i-1} \in p^{K-\tau(j)}R$, we see that
this is equivalent to 
$$
\big(-K + \tau(i) - \tau(2)\big) +
\big(\tau(j)+\tau(i-j+1) - \tau(i)\big) + \big(K-\tau(j)\big) \geq 0
$$
when $j > 0$.
But this is clear since the left hand side is just $\tau(i-j+1) -\
\tau(2)$.
For $j=0$, we need $\big(-K + \tau(i) - \tau(2)\big) + \big(L +
\tau(i+1) - \tau(i) - \tau(2)\big) \geq 0$.
Since $L_{i-1} \geq K+1$, it suffices to show that $1 + \tau(i+1) \geq
2 \tau(2)$ and this too is clear.
So $z_i \in R$ and as $z_i \in (x^i,y^i)R [p^{-1}]  \cap R$, it
follows that $z_i \in \big((x^i,y^i):_R p^m\big)$ for some $m$. 
Of course, $N=m$ suffices by hypothesis and so (3) holds as desired.

Now let $Q_i^*=\big(x^i,y^i, (xy)^{i-1}z_1, (xy)^{i-2}z_2, \ldots,
xyz_{i-1}\big)R$.
Let $D_i$ be the least nonnegative integer such that $p^{D_i}z_i \in
Q_i^*$ (certainly $D_i \leq N$).
Next let $L_i=L_{i-1}+ D_i$ and let $a_{ji0} = \Big(\prod_{m=0}^{j-1}
\frac{p^{L_i}-m}{p^{L_{i-1}}-m} \Big) a_{j,i-1}$.
By (2.5), we have 
$$
\sum_{j=0}^k\binom{p^{L_i}-j}{k-j}a_{ji0}
z^{k-j}x^j \in y^k R[p^{-1}]
$$ 
for $k < i$.
We also see that 
$$
\sum_{j=0}^{i-1}\binom{p^{L_i}-j}{i-j} a_{ji0}z^{i-j}x^j =
\big(p^{D_i}u\big)\big(p^{E_i}z_i\big) 
$$
where $u$ is the unit $\Big(\prod_{m=1}^{i-1}
\frac{p^{L_i}-m}{p^{L_{i-1}}-m} \Big)$.
Since $p^{D_i}z_i \in Q_i^*$, we have a relation 
$$
p^{D_i}z_i + cx^i + \sum_{n=1}^{i-1} c_n (xy)^{i-n}z_n = by^i
$$ 
and so 
$$
p^{E_i}u\bigg(p^{D_i}z_i + cx^i + \sum_{n=1}^{i-1}c_n(xy)^{i-n}z_n\bigg) \
\in\ y^iR[p^{-1}]\ .
$$  
Letting $a_{ii0} = p^{E_i}uc$, we have 
$$
\sum_{j=0}^i \binom{p^{L_i}-j}{i-j}a_{ji0}z^{i-j}x^j =
p^{E_i}u\big(p^{D_i}z_i+ cx^i\big) \ .
$$  
For each $n=1,\ldots,i-1$, we intend to find
$a_{1in},\ldots,a_{i-1,in}$ such that 
$$
\sum_{j=1}^{i-1} \binom{p^{L_i}-j}{i-j}a_{jin}z^{i-j}x^j =
p^{E_i}uc_n(xy)^{i-n}z_n \ .
$$  
Then we set $a_{ji} = \sum_{n=0}^{i-1}a_{jin}$ for $0 < j < i$ and
$a_{ii} = a_{ii0}$.
Clearly (1) will be proved if we show ($1'$) $a_{jin}$ is in $p^{K -
  \tau(j)}R$ for every $j$, $n$.
To prove (2), we note that the $i$th condition follows from the
definition:
\begin{equation*}
\begin{split}
\sum_{j=0}^i \binom{p^{L_i}-j}{i-j} a_{ji} z^{i-j}x^j & = 
\sum_{j=0}^i \binom{p^{L_i}-j}{i-j} a_{ji0}z^{i-j}x^j + 
\sum_{n=1}^{i-1}\sum_{j=1}^{i-1}\binom{p^{L_i}-j}{i-j}
a_{jin}z^{i-j}x^j \\
 & = p^{E_i}u (p^{D_i}z_i + cx^i) + \sum_{n=1}^{i-1} p^{E_i} u
 c_n(xy)^{i-n} z_n \in y^i  R[p^{-1}]\ .
\end{split}
\end{equation*}
Thus (2) will follow if we show ($2'$) $\sum_{j=0}^k \binom{p^{L_i}
  -j}{k-j} a_{jin} z^{k-j}x^j \in y^k R[p^{-1}]$ for $0 \leq n <
i$, $0 < k <i$.

We shall define the set $\{a_{jin}\}$ and prove ($1'$) and ($2'$) using
three cases: $n=0$, $n=1$, $n>1$.
For $n=0$, $\{a_{jin}\}$ is already defined and ($2'$) was observed
above.
For $j < i$, $a_{ji0} \in a_{j, i-1}R \subseteq p^{K-\tau(j)} R$ as
desired and finally $a_{ii0} = p^{E_i} uc \in p^{K-\tau(i) + \tau(2)}R
\subset p^{K-\tau(i)} R$, giving ($1'$).
For $n=1$, we get  
$$
\sum_{j=1}^{i-1} \binom{p^{L_i}-j}{i-j}a_{jin}z^{i-j}x^j =
p^{E_i}uc_n(xy)^{i-1} z
$$
by choosing  
$$
a_{jin}=\frac{p^{E_i}uc_ny^{i-1}}{\binom{p^{L_i}-j}{i-j}}
$$
when $j=i-1$ and $a_{jin}=0$ otherwise.  
That 
$$
\sum_{j=1}^k \binom{p^{L_i}-j}{k-j}a_{jin}z^{k-j}x^j \ \in \
y^kR[p^{-1}] 
$$ 
for every $k<i$ is trivial since each $a_{jin} \in
y^{i-1}R[p^{-1}]$.  
To see that $a_{jin} \in p^{K-\tau(j)}R$, we note that $a_{jin}=0$
unless $j=i-1$.   
In the latter case, we need only show 
$$
\frac{p^{E_i}}{\binom{p^{L_i}-(i-1)}{1}} \ \in\ p^{K-\tau(i-1)}R\ .
$$  
This requires only  
$
\big(K-\tau(i)+\tau(2)\big) - \big(\tau(i-1)+ \tau(2)-\tau(i)\big)
\geq K - \tau(i-1)
$
and in fact equality is clear.

Now we fix $n>1$.  
To get 
$$
\sum_{j=1}^{i-1}\binom{p^{L_i}-j}{i-j}a_{jin}z^{i-j}x^j =
p^{E_i}uc_n(xy)^{i-n}z_n\ ,
$$ 
we recall that
$$
z_n=p^{-E_n}\sum_{j=0}^{n-1}\binom{p^{L_{n-1}}-j}{n-j}
a_{j,n-1}z^{n-j}x^j \ .
$$  
We can obtain the desired inequality provided 
$$
\binom{p^{L_i}-j}{i-j}a_{jin} = p^{E_i}uc_ny^{i-n}p^{-E_n}
\binom{p^{L_{n-1}}-j-n+i}{i-j} a_{j+n-i,n-1}
$$ 
for $j=i-n, \ldots,
i-1$  and $a_{jin}=0$
otherwise.  
So for $j=i-n, \ldots, i-1$, 
\begin{equation*}
\begin{split}
a_{jin} & =
\frac{p^{E_i}uc_ny^{i-n}p^{-E_n}\binom{p^{L_{n-1}}-j-n+i}{i-j}
  a_{j+n-i,n-1} }{\binom{p^{L_i}-j}{i-j}} \\
 & =(p^{E_i-E_n}uc_n)
\frac{y^{i-n}\binom{p^{L_{n-1}}-j+(i-n)}{i-j}a_{j-(i-n),n-1}}
{\binom{p^{L_i}-j}{i-j}}  
\end{split}
\end{equation*} 
Applying (2.6), we see that ($2'$) holds (use $i-n$ for $d$, $L_i$ for
$M$, $L_{n-1}$ for $L$, and $a_{j,n-1}$ for $a_j$).

Finally, to prove that $a_{jin} \in p^{K-\tau(j)}R$, it suffices to
show that 
$$
\frac{p^{E_i-E_n}\binom{p^{L_{n-1}}-j-n+i}{i-j}
  a_{j+n-i,n-1}}{\binom{p^{L_i}-j}{i-j}}\ \in \  p^{K-\tau(j)}R\ . 
$$  
For $j+n-i>0$, we apply (\ref{lem:1.6}) to see that this is equivalent
to 
\begin{multline*}
\big(K - \tau(i) + \tau(2)\big) - \big(K - \tau(n)+ \tau(2)\big) +
\tau(j+n-i) + \tau(i-j+1) - \tau(n) \\
+ \big(K-\tau(j+n-i)\big) - \big(\tau(j)+
\tau(i-j+1)-\tau(i)\big)\geq K -\tau(j)\ . 
\end{multline*}
However, the two sides of this expression are clearly equal.  
For $j+n-i=0$, we must show 
$$
\frac{p^{E_i-E_n}\binom{p^{L_{n-1}}}{i-j}}{\binom{p^{L_i}-j}{i-j}} \ \in
\ p^{K-\tau(j)}R \ .
$$  
Using (\ref{lem:1.6}), it suffices to show 
\begin{multline*}
\big(K - \tau(i) + \tau(2)\big) - \big(K - \tau(n)+ \tau(2)\big) + 
\big(L_{n-1} - \tau(i-j) + \tau(i-j+1) - \tau(2) \big) \\- \big(\tau(j)
+ \tau(i-j+1) - \tau(i)\big) \geq K -\tau(j)\ .
\end{multline*}
This is equivalent to $\tau(n) + L_{n-1} - \tau(2) -\tau(i-j) \geq K$.   
Since $i-j=n$, this is equivalent to $L_{n-1}\geq K + \tau(2)$.   
Since $\tau(2) \leq 1$, and $L_{n-1} \geq K+1$, this is clear.
Thus ($1'$) holds and we may satisfactorily perform the inductive step.

Finally, by (2.4) we see that for $i \gg 0$, $D_i = 0$ and so $L_i$ is
eventually constant.
Thus we eventually reach the $i=p^{L_i}$ step and the proof is
complete.
\end{proof}

\begin{rmk} 
If $N=1$ and $H_1(p,x,y; R)$ is cyclic, we may draw the stronger
conclusion that $z$ is in the plus closure of $(x,y)R$.  
The proof is a much simplified version of (2.3).  
We let $K=0$, $L=1$, and dispense with the $\tau$ function entirely.  
Everything works with far fewer steps.  
We do not know whether the plus closure captures the colon in the
three dimensional case in general.     
\end{rmk}

Now we state our corollaries.

\begin{cor} 
For the class of equidimensional three-dimensional excellent semilocal
domains of mixed characteristic, the full extended plus closure has
the colon-capturing property. 
\end{cor}

\begin{proof} 
Let $R$ be any equidimensional three-dimensional excellent semilocal
domain of mixed characteristic.  
We may adjoin a $(p-1)^{\text{st}}$ root of $p$ to $R$ if $R$ does not
already contain one.  
We also take the integral closure.  
Our hypotheses are preserved and it suffices to demonstrate
colon--capturing for the new domain.  
Choose $x,y$ such that $p,x,y$ is a complete system of
parameters.  
By (1.1), replacing $x,y$ by powers if necessary, we may
assume that $p^N,x,y$ kill all first Koszul homology modules.  
Then we can apply (2.7) to see that $\big((x,y):_R p^N\big)
\subset (x,y)^{\epf}$.  
The result now follows from (1.4).
\end{proof}

\begin{cor} \label{cor:2.7}
Let $R$ be a three-dimensional regular local ring of mixed
characteristic and let $I$ be an ideal of $R$. 
Then $I= I^+ = I^{\epf} = I^{\rf}$.
\end{cor}

\begin{proof} 
It suffices to show  $I=I^{\rf}$ since the other two closures
are trapped in between.  
If an element is in the full rank one closure of $I$, this remains
true after completing and so it suffices to prove the corollary for
complete rings.    
Now we combine (1.2) and (2.9).
\end{proof}

On behalf of the reader, I'd like to extend thanks to Craig Huneke,
Paul Roberts, and the referee for various comments which led to a
clearer presentation.
In particular, it was Paul Roberts \cite{R2} who brought to my
attention why (2.5) and (2.6) work, the connection with derivatives
and integrals.

\begin{center}
\bigskip{\it Department of Mathematics, University of Texas at Austin, Austin
  TX 78712} \\
E-mail: heitmann@math.utexas.edu
\end{center}

\end{document}